\documentstyle[12pt]{article}
\newtheorem{Theorem}{Theorem}[section]
\newtheorem{Lemma}[Theorem]{Lemma}
\newtheorem{Corollary}[Theorem]{Corollary}

\newtheorem{Definition}[Theorem]{Definition}
\newtheorem{Fact}[Theorem]{Fact}
\newtheorem{Example}[Theorem]{Example}

\def\Z{\mbox{\boldmath$Z$}}
\def\Q{\mbox{\boldmath$Q$}}
\def\x{\bar{x}}
\def\y{\bar{y}}
\def\m{\bar{m}}
\def\a{\bar{a}}
\def\b{\bar{b}}
\def\Av{{\rm Av}}
\def\Th{{\rm Th}}
\def\lor{\;\vee\;}
\def\land{\;\wedge\;}
\def\cnc{\hat{\ }}
\def\cset#1#2{\{#1 : #2\}}
\def\seq#1#2{\langle#1 : #2\rangle}

\newenvironment{Proof}{\noindent {\bf Proof}\quad}%
{\quad\vrule height7pt depth0pt width4pt \bigskip}

\title{Model Companions of Theories
with an Automorphism
\footnotetext{The work has been done while the author was
visiting MSRI. The author would like to thank MSRI for its
hospitality.}
}
\author{Kikyo, Hirotaka}
\date{June 29, 1998}

\begin{document}

\maketitle

\begin{abstract}
For a theory $T$ in $L$, $T_\sigma$ is the theory of
the models of $T$ with an automorphism $\sigma$.
If $T$ is an unstable model complete theory without the independence
property, then $T_\sigma$ has no model companion.
If $T$ is an unstable model complete theory and $T_\sigma$ has
the amalgamation property, then $T_\sigma$ has no model companion.
If $T$ is model complete and has the fcp, then 
$T_\sigma$ has no model completion.
\end{abstract}

\section{Introduction}

In \cite{ChP}, Chatzidakis and Pillay studied 
the theory $T_P$, the theory of models of $T$ with
a distinguished new unary predicate $P$,
and $T_\sigma$, the theory of models of $T$
with a distinguished new automorphism $\sigma$.
They showed that if $T$ is stable, 
$T = T^{eq}$, and admits QE, then
the model companions of $T_P$ and $T_\sigma$ are simple.
The typical case is that the model companion of
$ACF_\sigma$ exists and called $ACFA$,
where $ACF$ is the theory of algebraically closed field.
They also showed that
the model companion of $T_P$ exists
if and only if $T$ admits the elimination of $\exists^\infty$,
but they did not show any conditions on the existence
of the model companion of $T_\sigma$.

Khudaibergenov proved that if $T$ is stable and 
has the fcp then $T_\sigma$ has no model companion,
by reducing to 
the case that $T$ does not admit the elimination of $\exists^\infty$.
In this case the argument is similar to that for $T_P$.

In this paper, we show that $T_\sigma$ has no model companion
if $T$ is unstable without the independence property,
or $T$ is unstable and $T_\sigma$ has the amalgamation property.
If $T$ is stable then $T_\sigma$ has the amalgamation property,
or more generally it has PAPA \cite{Lascar1, Lascar2, ChP} in $T^{eq}$.
So, if $T$ is model complete and has the fcp, then
$T_\sigma$ has no model completion.
In particular,
if $T=T^{eq}$, admits the elimination of quantifiers, and has PAPA,
then $T_\sigma$ has no model companion if $T$ has the fcp.

Here is a history of the results.
Tsuboi tried to extend the result of Chatzidakis and Pillay
for $T_\sigma$ to the case $T$ is simple, and investigated
the theory of the random graph with an automorphism.
He tried to give an axiomatization of the model companion
of the theory, and found some difficulty when $\sigma$ has
a fixed point of certain kind.
Looking at this difficulty, the author found that the theory of
the random graph with an automorphism does not have
the model companion. When the author told Pillay about this,
he suggested that the fcp might be the reason.
The author checked this property for the theory of
the dense linear order without endpoints, the o-minimal theories,
and the theory of the atomless Boolean algebra. It turned out that
there are no model companions of $T_\sigma$
for these theories.
With some more work, 
we could generalize them to the theorems in this paper.

We could not get rid of the assumption that $T_\sigma$ has
the amalgamation property in general.
For our argument, the formulas with the independence property
are harmful. The random graph has the independence property
and the unstability comes only from this, but in this case
the argument for the existence of an appropriate fixed point
works because of its randomness. 
Even if we assume that $T$ is simple
or supersimple, it seems that the machineries developed
for the simple theories don't give good information.

There are some theories $T$ such that
$T_\sigma$ has no model companion but
has independence property and 
does not have the amalgamation property.
ACFA is such an example.
It seems that our argument does not work
in this example.

The author would like to thank Akito Tsuboi
for the valuable discussions about the random graph,
and Zo\'e Chatzidakis and Ehud Hrushovski for the
information about ACFA and pseudo-finite fields.
The author also would like to thank
Daniel Lascar and Anand Pillay 
for the valuable discussions and comments.

\section{Preliminaries}

First we recall some basic definitions and facts.
If $M \subset N$ are models of a model complete theory
then $M$ is an elementary substructure of $N$.
So the truth of any definable relations with
parameters from $M$ are preserved upwards and downwards. In fact,
the following are equivalent:
\begin{enumerate}
\renewcommand{\labelenumi}{(\arabic{enumi})}
\item $T$ is model complete;
\item Every formula is equivalent to an existential formula modulo $T$;
\item Every formula is equivalent to a universal formula modulo $T$.
\end{enumerate}
Note also that a model complete theory can be axiomatized by
universal-existential sentences (a universal-existential theory).
If a theory $T$ admits the elimination of quantifiers then $T$ is
model complete.

\begin{Definition}{\rm
Let $T$ and $T'$ be two theories (not necessarily complete) 
in the same language.
$T'$ is a {\em model companion} of $T$ if
$T'$ is model complete, every model of $T$ can be extended to
a model of $T'$ and every model of $T'$ can be extended to
a model of $T$.
$T'$ is a {\em model completion} of $T$ if
$T'$ is a model companion of $T$, and
if $M$ is a model of $T$ and $M_1$, $M_2$ are models of $T'$
extending $M$ then $M_1$ and $M_2$ are elementarily
equivalent over $M$.
If $T$ has a model completion then $T$ must have the
amalgamation property and if $T$ has the amalgamation property
then the model companion of $T$ is a model completion of $T$.
}\end{Definition}

Recall that if $T$ is a universal-existential theory
and $T'$ a model companion of $T$, then
$M$ is a model of $T'$ if and only if $M$ is an existentially 
closed model of $T$.

We need some more notation and definition for our arguments.
$\langle\cdots\rangle$
denotes a sequence,
and $I\cnc I'$ denotes the
concatenation of sequences $I$ and $I'$.

\begin{Definition}{\rm
Let $I = \seq{\a_i}{i<\omega}$ be a sequence of tuples,
$L$ a language and
$B$ a set of parameters.
The {\em average $L$-type} of $I$ over $B$ is the partial type
$\Av_L(I/B)$ defined by
\[
\Av_L(I/B) = \cset{\varphi(\x,\b)}
{\models \varphi(\a_i,\b)\mbox{ for all but finitely many }i<\omega}.
\]
}\end{Definition}

\begin{Definition}{\rm 
Let $M$ be an $L$-structure, $\sigma$ an automorphism of $M$ and 
$p(\x)$ an $L$-type over $M$.
$\sigma p(\x)$ is the type
\[\cset{\varphi(\x,\sigma(\m))}{\varphi(\x,\m)\in p(\x),\; 
\varphi(\x,\y) \in L}.
\]
}\end{Definition}

\section{Unstable models with an automorphism}

For the rest of the paper,
unless stated otherwise, we assume that $T$ is a complete 
and model complete theory in the language $L$.
We use $\sigma$ to denote the unary function symbol representing
an automorphism of a model of $T$.
$T_\sigma$ is
the union of $T$ and the set of axioms 
stating that $\sigma$ is an $L$-automorphism.
Note that $T_\sigma$ is a universal-existential
theory in $L(\sigma) = L\cup\{\sigma\}$.

% which is 
%and $\x < \y$ is a definable relation in $L$ which gives
%an infinite order in some model of $T$, i.e.,
%for some $\{\a_i\,|\,i<\omega\}$, $\a_i < \a_j$ if
%and only if $i < j$. We can assume that $<$ is antisymmetric.

\begin{Theorem}\label{th:NIP}
If $T$ is an unstable model complete theory without the
independence property, then $T_\sigma$ has no model companion.
\end{Theorem}

\begin{Proof}
Since $T$ is unstable, there is an 
antisymmetric $L$-definable relation
$\x < \y$ such that for some sequence $\cset{\a_i}{i<\omega}$,
$\a_i < \a_j$ if and only if $i<j$.

  Suppose that $T_\sigma$ has the model companion $T_A$.
Consider the following theory $T_1$ in the language $L(\sigma,\a,\b)$
where $\a$ and $\b$ are sequences of new constants:
\begin{enumerate}
\renewcommand{\labelenumi}{(\arabic{enumi})}
\item The axioms in $T_A$;
\item $\seq{\sigma^i(\a)}{i<\omega}\cnc\langle\b\rangle$ 
is $L$-indiscernible;
\item $\a < \sigma(\a)$.
\end{enumerate}

We claim first that $T_1$ is consistent.
First, we work in a big model of $T$.
Let $\seq{\a_i}{i<\omega}\cnc\langle\b\rangle$ be an $L$-indiscernible
sequence such that $\a_0 < \a_1$.
We can choose a model $M$ of $T$ and 
an automorphism $\sigma$ of $M$ such that 
$M$ contains this sequence and
$\sigma(\a_i) = \a_{i+1}$ for $i < \omega$. 
$M$ satisfies (2) and (3) with $\a = \a_0$.
Since $(M,\sigma)$ is a model of $T_\sigma$,
it can be extended to a model $(M',\sigma')$ of $T_A$.
Since $M'$ is an $L$-elementary extension of $M$,
$M'$ satisfies (2) and (3) also.

Now we claim that 
$T_1 \vdash \exists \x\,[\sigma(\x)=\x \land \a < \x < \b]$.

Let $(M,\sigma)$ be any model of $T_1$.
Then $I = \seq{\sigma^i(\a)}{i<\omega}$ is an $L$-indiscernible
sequence in $M$. Let $p(\x) = \Av_L(I/M)$. Since $T$ does not have
the independence property, $p(\x)$ is a complete $L$-type over $M$.
By the definition of the average type, 
the formula $\a < \x < \b$ belongs to $p(\x)$.
Also, it is easy to check that $\sigma(p) = p$.
So, there is an elementary extension $M'$ of $M$ and an
automorphism $\sigma'$ of $M'$ such that
$\sigma'$ extends $\sigma$ and fixes some realization of $p(\x)$ in $M'$.
In particular, 
\[
(M',\sigma') \models \exists \x\,[\sigma(\x)=\x \land \a < \x < \b].
\]
Since $(M,\sigma)$ is existentially closed among the models of
$T_\sigma$ and the relation $<$ can be written in an existential form
modulo $T$,
\[
(M,\sigma) \models \exists \x\,[\sigma(\x)=\x \land \a < \x < \b].
\]
Thus we have the claim.

By compactness,
$T'_1 \vdash \exists \x\,[\sigma(\x)=\x \land \a < \x < \b]$ for some
finite subset $T'_1$ of $T_1$.
In particular, the axiom (2) can be replaced by
the indiscernibility of 
$\seq{\sigma^i(\a)}{i<N}\cnc\langle\b\rangle$ for some natural
number $N$. 
Now choose a model $(M,\sigma,\a,\b_0)$ of $T_1$.
If we change the interpretation of $\b$ to
$\sigma^k(\a)$ for some $k > N$, still it is a model of $T'_1$.
We have
\[
(M,\sigma,\a,\sigma^k(\a))\models \exists \x\,[\sigma(\x)=\x \land \a < \x < \sigma^k(\a)],
\]
but this cannot happen because $\sigma$ 
is an $L$-automorphism and must preserve the relation $<$.
\end{Proof}

The proof of the theorem above is a generalization of
the direct proof of the fact that the following examples have no
model companion.

\begin{Example}{\rm 
The model companions of $T_\sigma$ do not exist for
the following theories:
The theory $T$ of dense linear order, or
more generally,
model complete o-minimal theories;
The theory of $p$-adic number field in Macintyre's language.
}\end{Example}

There are many theories which have the independence property
but the corresponding $T_\sigma$ has no model companion.

\begin{Theorem}\label{th:AP}
If $T$ is an unstable model complete theory and $T_\sigma$ has
the amalgamation property, then $T_\sigma$ has no model companion.
\end{Theorem}

\begin{Proof}
Let $T_S$ be a Skolemization of $T$ in the language $L_S$.

First, we work in a big model of $T_S$.
Since $T$ is unstable, there is an 
antisymmetric $L$-definable relation
$\x < \y$ such that for some sequence $\cset{\a_i}{i<\omega}$,
$a_i < a_j$ if and only if $i<j$.
Choose an $L_S$-indiscernible sequence
$\seq{\a_i}{i\in \Z}\cnc\seq{\b_i}{i\in \Z}$
such that $\a_0 < \a_1$.
Let $M_0$ be the Skolem hull of this sequence.
Then the map $\sigma_0$ given by 
$\sigma_0(\a_i)=\a_{i+1}$ and $\sigma_0(\b_i)=\b_{i+1}$
extends to an $L_S$-automorphism of $M_0$.

The reduct of $(M_0,\sigma_0)$ to $L(\sigma)$ is a model of $T_\sigma$.

Now suppose that $T_\sigma$ has the model companion $T_A$.
Then $(M_0,\sigma_0)$ has an extension
$(M_1,\sigma_1)$ which is a model of $T_A$.

Now consider the following theory $T_1$ in
the language $L_S(\sigma,P,\a,\b)$ where $P$ is a new unary
predicate symbol, and $\a$, $\b$ are sequences of new constant
symbols:
\begin{enumerate}
\renewcommand{\labelenumi}{(\arabic{enumi})}
\item The axioms in $T_A$;
\item Every functions in $L_S$ are total on $P$ and
$P$ is closed under every functions in $L_S$;
\item $P$ is a model of $T_S$;
\item $\sigma$ is an $L_S$-automorphism on $P$;
\item \label{ax:indis} 
$\seq{\sigma^i(\a)}{i\in\Z}\cnc\seq{\sigma^i(\b)}{i\in\Z}$ 
is $L_S$-indiscernible in $P$;
\item $\a < \sigma(\a)$.
\end{enumerate}
Note that the Skolem functions in $L_S\setminus L$ are considered here
as partial functions.
If we interpret $P$ by $M_0$ in $M_1$, we get a model of $T_1$.
Thus $T_1$ is consistent.

Now we claim that 
$T_1 \vdash \exists \x\,[\sigma(\x)=\x \land \a < \x < \b]$.

Let $(M_1,\sigma_1,\a,\b)$ be any model of $T_1$.
Let $M_0$ be the Skolem hull of the sequence
$I = \seq{\sigma^i(\a)}{i\in\Z}\cnc\seq{\sigma^i(\b)}{i\in\Z}$
in $P(M_1)$. Note that $M_0$ is an elementary $L_S$-substructure of
$P(M_1)$. 
In particular, $M_0$ is a model of $T$.

Since $\sigma_1$ is an $L_S$-automorphism on $P(M_1)$,
$M_0$ is closed under $\sigma_1$ and $\sigma_0 = \sigma_1|M_0$ is an 
$L_S$-automorphism of $M_0$.
Hence, $(M_0,\sigma_0)$ is a model of $T_\sigma$.

Now consider the average $L$-type $p(\x)$ of $\seq{\sigma_0^i(\a)}{i < \omega}$
over $M_0$. 
%i.e.,
%for $\varphi$ in $L$ and $\m$ in $M_0$,
%$\varphi(\x,\m)$ belongs to $p(\x)$ if and only if
%for sufficiently large $n$, $M_0 \models \varphi(\a_n,\m)$.
Let $\varphi(\x,\m)$ be an $L$-formula with parameters $\m$ in $M_0$.
Since $M_0$ is the Skolem hull of $I$, 
$\varphi(\x,\m)$ can be written as
a formula in $L_S$ with finite number of parameters from $I$,
say $\sigma_0^i(\a)$ and $\sigma_0^i(\b)$ for $|i| < k$.
By the $L_S$-indiscernibility of $I$,
the truth value of $\varphi(\sigma_0^n(\a),\m)$ are the same for $n > k$.
So, one of $\varphi(\x,\m)$ or $\neg\varphi(\x,\m)$ belongs to $p(\x)$.
Hence, $p(\x)$ is complete over $M_0$.
Note that the formula $\a < \x < \b$ belongs to $p(\x)$.
Also, it is easy to check that $\sigma_0(p) = p$.

Now we can find an $L$-elementary extension $M_2$ of $M_0$
and an $L$-auto\-mor\-phism $\sigma_2$ of $M_2$ such that
$\sigma_2$ extends $\sigma_0$ and fixes some realization 
of $p(\x)$ in $M_2$. In particular,
$(M_2,\sigma_2)$ is a model of $T_\sigma$ and
\[
(M_2,\sigma_2) \models \exists \x\,[\sigma(\x)=\x \land \a < \x <\b].
\]
By the amalgamation property of $T_\sigma$,
there is an extension $(M_3,\sigma_3)$ of $(M_1,\sigma_1)$ such that
$(M_3,\sigma_3)$ is a model of $T_\sigma$ and 
$(M_2,\sigma_2)$ embeds in $(M_3,\sigma_3)$ over $(M_0,\sigma_0)$.
Then we have
\[
(M_3,\sigma_3) \models \exists \x\,[\sigma(\x)=\x \land \a < \x <\b].
\]
Since $T$ is model complete, $<$ can be written in an existential form
modulo $T$. 
As $(M_1,\sigma_1)$ is a model of $T_A$, it is existentially closed
among the models of $T_\sigma$. Hence,
\[
(M_1,\sigma_1) \models \exists \x\,[\sigma(\x)=\x \land \a < \x <\b].
\]
Thus we have the claim.

By compactness, only a finite part of $T_1$ is necessary to deduce
the formula $\exists \x\,[\sigma(\x)=\x \land \a < \x <\b]$.
In particular, the axiom (\ref{ax:indis}) can be replaced by the axiom
\[
\seq{\sigma^i(\a)}{|i| < N}\cnc\seq{\sigma^i(\b)}{|i| < N}
\mbox{ is $L_S$-indiscernible in $P$};
\]
\noindent for some natural number $N$ and still deducing the above formula.
We can satisfy these axioms
by taking a model $(M,\sigma,\a,\b_0)$ of $T_1$ and
changing the interpretation of $\b$ to $\sigma^k(\a)$ for some $k > 2N$.
Then we have 
\[
(M,\sigma,\a,\sigma^k(\a))\models 
\exists \x\,[\sigma(\x)=\x \land \a < \x <\sigma^k(\a)],
\]
but this cannot happen because $\sigma$ 
is an $L$-automorphism and must preserve the relation $<$. 
\end{Proof}

With the following fact proved by Kudaibergenov,
we can conclude that if $T$ has the fcp and model complete
then $T_\sigma$ has no model completion.

\begin{Definition}{\rm
For a complete theory $T$, we say that
{\em $T$ does not admit the elimination of $\exists^\infty$}
if there are a formula
$\varphi(\x,\y)$ and a set of parameters 
$\cset{\a_n}{n<\omega}$ such that
for each $n <\omega$, the number of the realizations 
of $\varphi(\x,\a_n)$ is finite but greater than $n$.
}\end{Definition}

\begin{Fact}[Kudaibergenov]
\begin{enumerate}
\renewcommand{\labelenumi}{{\rm (\arabic{enumi})}}
\item \label{facta} If $T$ is model complete and does not admit 
the elimination of $\exists^\infty$,
then $T_\sigma$ has no model companion.
\item If $T$ is stable, model complete and has the fcp,
then $T_\sigma$ has no model companion.
\end{enumerate}
\end{Fact}

\begin{Proof} (1) Suppose $T$ does not admit the
elimination of $\exists^\infty$ and $T_\sigma$ has the
model companion $T_A$.
Let $\varphi(\x,\y)$ and a set of parameters 
$\cset{\a_n}{n<\omega}$ be such that
for each $n <\omega$, the number of the realizations 
of $\varphi(\x,\a_n)$ is finite but greater than $n$.
Let $M_0$ be a model of $T$ containing $a_n$'s.
$(M_0,id)$ where $id$ is the identity map on $M$ is a model
of $T_\sigma$, so there is a model $(M,\sigma)$ of $T_A$
extending $(M_0,id)$. Since $M_0 \prec M$ and
the set of realizations of each $\varphi(\x,\a_n)$ in $M_0$
is finite, it is the same in $M$. So, $\sigma$ fixes
every realization of each $\varphi(\x,\a_n)$.
Let ${\cal U}$ be a non-principal ultrafilter on $\omega$.
Make an ultrapower 
$(\tilde{M},\tilde{\sigma}) = (M,\sigma)^\omega/{\cal U}$,
and let $\b = \seq{\a_i}{i<\omega}/{\cal U}$.
Since each $\sigma$ is the identity on each $\varphi(\x,\a_i)$,
$\tilde{\sigma}$ is also the identity on $\varphi(\x,\b)$.
On the other hand, 
the number of realizations of $\varphi(\x,\b)$ in $\tilde{M}$
is infinite. 
But it is impossible to have a model of $T_A$ such that
the automorphism fixes pointwise an infinite definable set.

(2) Suppose $T$ is stable, model complete, and has the fcp.
Then $T^{eq}$ is also model complete and 
does not admit the elimination of $\exists^\infty$
by the fcp theorem \cite{Shelah}.
Suppose in addition that $T_\sigma$ has the model companion
$T_A$. Now it is routine to check that
$T_A\cup (T^{eq})_\sigma$ is the model companion of
$(T^{eq})_\sigma$. This contradicts (1).
\end{Proof}

\begin{Corollary}
Suppose $T$ is model complete and not necessarily complete.
If $T$ has some model with the fcp, then $T$ has no model completion.
In particular, if $T= T^{eq}$, admits the quantifier elimination,
has the PAPA \cite{Lascar1, Lascar2, ChP}
and the fcp, then $T_\sigma$ has no model companion.
\end{Corollary}

\begin{Example}{\rm The theories of
the random graph and the atomless Boolean algebra
have the PAPA. More generally, unstable countably categorical
theories obtained by the Fra\"{\i}ss\'e's method 
from a class of finite structures which admits
the free amalgamation have the PAPA \cite{Lascar2}. 
The model companions of $T_\sigma$ do not exist for
these theories.
}\end{Example}

\section{Some remarks}

Hrushovski proved that $ACFA_\sigma$ does not have 
neither the model companion nor the amalgamation property.
Chatzidakis proved that $Psf_\sigma$,
where $Psf$ is the theory of pseudo-finite fields,
does not have the model companion by a similar argument to Hrushovski's.
Their arguments are very specific to the field theory.
It seems that our method does not work in these theories.

Here is another example of a theory $T$ such that
$T_\sigma$ does not have the amalgamation property
due to M.~Ziegler \cite{Lascar1}.

\begin{Example}{\rm The language has two predicate symbols
$R(x,y,z)$ and $U(x)$. Let $(M, <)$ be a dense linear
order without endpoints and interpret $R$ and $U$ on $M$ so that
$R(a,b,c)$ if and only if $a < b < c$ and $U(M)$
is dense and co-dense in $M$.
Then $\Th(M;R,U)$ admits the quantifier
elimination, but does not have 
neither the amalgamation property
nor the independence property.

Hence, the theorems \ref{th:NIP} and \ref{th:AP} are independent.
}\end{Example}

It is possible to expand the Ziegler's example to have the independence 
property without the $T_\sigma$ having the amalgamation property.
To show this, we give a lemma concerning the types
fixed by an automorphism. The proof is similar to that of
the existence of a nonforking extension of a type in \cite{Shelah}.

\begin{Definition}{\rm
Let $M$ be a model and $\sigma$ an automorphism of $M$.
A formula $\varphi(\x,\m)$, where $\varphi(\x,\y)$ is parameter-free
and $\m \in M$, is said to be {\em $\sigma$-consistent}
if $\cset{\varphi(\x,\sigma^n(\m))}{n<\omega}$ is consistent.
}\end{Definition}

\begin{Lemma}Let $M$ be a model, $\sigma$ an automorphism of $M$,
and $p(\x)$ a partial type over $M$. The following are
equivalent.
\begin{enumerate}
\renewcommand{\labelenumi}{{\rm (\arabic{enumi})}}
\item There is a complete type $q(\x)$ over $M$ such that
$p(\x)\subset q(\x)$ and $\sigma q(\x) = q(\x)$;
\item For any formulas $\varphi_1(\x,\m_1)$, $\ldots$,
$\varphi_k(\x,\m_k)$ over $M$,
\[
p(\x) \vdash \varphi_1(\x,\m_1)\lor\cdots\lor\varphi_k(\x,\m_k)
\]
implies that $\varphi_i(\x,\m_i)$ is $\sigma$-consistent for some $i$
between $1$ and $k$. 
\end{enumerate}
\end{Lemma}

\begin{Proof} (1) $\Rightarrow$ (2). 
Let $q(\x)$ be a complete type over $M$ such that
$p(\x)\subset q(\x)$ and $\sigma q(\x) = q(\x)$.
Suppose
\[
p(\x) \vdash \varphi_1(\x,\m_1)\lor\cdots\lor\varphi_k(\x,\m_k).
\]
Then one of $\varphi_i(\x,\m_i)$'s is in $q(\x)$.
It is $\sigma$-consistent.

(2) $\Rightarrow$ (1). 
Let 
\[
\Sigma(\x) = \cset{\neg \varphi(\x,\m)}
{\varphi(\x,\m)\mbox{\ is not $\sigma$-consistent,}
\;\varphi(\x,\y)\in L,\;\m\in M}.
\]
By compactness and (2), $p(\x) \cup \Sigma(\x)$ is consistent.
Let $q(\x)$ be any completion of $p(\x) \cup \Sigma(\x)$ over $M$.

Suppose $\sigma q(\x) \neq q(\x)$. Then $\varphi(\x,\m) \in q(\x)$
and $\neg\varphi(\x,\sigma(\m)) \in q(\x)$, thus
$\varphi(\x,\m)\land \neg\varphi(\x,\sigma(\m)) \in q(\x)$
for some $\varphi(\x,\m)$.
On the other hand,
$\varphi(\x,\m)\land \neg\varphi(\x,\sigma(\m))$ is
an $L(M)$-formula but is not $\sigma$-consistent. So, its
negation belongs to $\Sigma(\x) \subset q(\x)$.
This is a contradiction. Hence $\sigma q(\x) = q(\x)$.
\end{Proof}

It seemed that this lemma
would be useful for our general argument to find a fixed point, 
but we could not use it.

\begin{Example}{\rm 
There are many expansions of the Ziegler's example
whose $T_\sigma$ do not have the amalgamation property.

Let $M = \Q\setminus \{0\}$ and 
consider $(M;R(x,y,z),U(x))$ where 
$R(a,b,c)$ if and only if $b$ is in between $a$ and $c$, 
and $U(M) = \cset{m/2^n}{m,n\in\Z}$. 

Let $\alpha$ be the automorphism of 
$M$ given by $\alpha(x) = -x$ for $x\in \Q$.
Put any relations on $\Q^+ = \cset{x\in\Q}{x>0}$.
Expand the each relation $P(x_1,\ldots,x_k)$ to $M$
in the minimum way so that
\begin{itemize}
\item $P(x_1,\ldots,x_k)$ if and only if  
$P(\alpha(x_1),\ldots,\alpha(x_k))$, and

\item $P(x_1,\ldots,x_i,\ldots,x_k)$ if and only if  
$P(x_1,\ldots,x_{i-1},\alpha(x_i),x_{i+1},\ldots,x_k)$ 
for each $i$.
\end{itemize}

$\alpha$ is an automorphism for the expanded structure.
Also,
each new relation $P$ restricted to $\Q^+$ does not change
by this modification. So, if you put a relation with
the independence property at the beginning, the property
is preserved.

Now consider the two partial types
\begin{eqnarray*}
p_1(x)&=&\cset{R(a,x,b)}{a < 0 < b}\cup \{U(x)\}
\mbox{ and}\\
p_2(x)&=&\cset{R(a,x,b)}{a < 0 < b}\cup \{\neg U(x)\}.
\end{eqnarray*}
By the lemma above, it is easy to check that
both $p_1(x)$ and $p_2(x)$ have some complete extensions
to $M$ which are fixed by $\alpha$.
We can extend $\alpha$ in two ways so that
it has some fixed point realizing each one of $p_1(x)$ and $p_2(x)$,
but it is impossible to have two distinct fixed point at the same time
extending $\alpha$.
}\end{Example}

If we do not put any new relations on some open interval containing 0 
then there will be no model companions for $T_\sigma$,
but we don't know if this is true in general.

\begin{tabbing}
Department of Mathematical Sciences\\
Tokai University\\
1117 Kitakaname, Hiratsuka 259-12\\
Japan\\
e-mail:{\tt kikyo@ss.u-tokai.ac.jp}
\end{tabbing}

\end{document}